\documentclass[jip]{degruyter-journal-a}      % J. Inverse Ill-Posed Probl.

\usepackage{amssymb}
\usepackage{amsmath}
\usepackage{amsfonts}
\usepackage{amsthm}
\usepackage{times}
\usepackage{comment}
\usepackage{color}
\usepackage{mathrsfs}
\usepackage{graphicx}
\usepackage{float}

\usepackage{tikz}
\usetikzlibrary{positioning}

\numberwithin{equation}{section}

\newtheorem{thm}{Theorem}[section]
\newtheorem{prop}[thm]{Proposition}
\newtheorem{lem}[thm]{Lemma}

\theoremstyle{definition}
\newtheorem{defn}[thm]{Definition}

\numberwithin{equation}{section}

\newcommand{\R}{\mathbb{R}}  % The real numbers.
\newcommand{\N}{\mathbb{N}}  % The Natural numbers.
\newcommand{\sw}{\mathcal{S}} % The Schwartz space
\newcommand{\domain}{[0,2\pi]\times \R} % domain

\title{Modified Radon transform inversion using moments}
\headlinetitle{Modified Radon transform inversion using moments}

\lastnameone{Choi}
\firstnameone{Hayoung}
\nameshortone{H.~Choi}
\addressone{School of Information Science and Technology, 
Shanghai Tech University, 
Shanghai 201210}
\countryone{China}
\emailone{hchoi@shanghaitech.edu.cn}

\lastnametwo{Ginting}
\firstnametwo{Victor}
\nameshorttwo{V.~Ginting}
\addresstwo{Department of Mathematics and Statistics, University of Wyoming, Laramie, WY 82071-3036}
\countrytwo{USA}
\emailtwo{vginting@uwyo.edu}

\lastnamethree{Jafari}
\firstnamethree{Farhad}
\nameshortthree{F.~Jafari}
\addressthree{Department of Mathematics and Statistics, University of Wyoming, Laramie, WY 82071-3036}
\countrythree{USA}
\emailthree{fjafari@uwyo.edu}

\lastnamefour{Mnatsakanov}
\firstnamefour{Robert}
\nameshortfour{R.~Mnatsakanov}
\addressfour{Department of Mathematics, West Virginia University
Morgantown, WV 26506}
\countryfour{USA}
\emailfour{robert.mnatsakanov@math.wvu.edu}

\abstract{Moment methods to reconstruct images from their Radon transforms are both natural and useful. They can be used to suppress noise or other spurious effects and can lead to highly efficient reconstructions from relatively few projections. We establish a modified Radon transform (MRT) via convolution with a mollifier and obtain its inversion formula. The relationship of the moments of the Radon transform and the moments of its modified Radon transform is derived and MRT data is used to provide a uniform approximation to the original density function. The reconstruction algorithm is implemented, and a simple density function is reconstructed from moments of its modified Radon transform. Numerical convergence of this reconstruction is shown to agree with the derived theoretical results.}



\begin{document}

\section{Introduction} \label{sec:intro}
%%%%%%%%%%%%%%%%%%%%%%%%%%%%%%%%%%%%%%%%
Radon transform of an integrable function over $\R^2 $ is the integral of that function over lines. A key application of Radon transform is tomography where the interior density of a 2-D object (e.g. slices of a 3-D object) is reconstructed from its Radon transform data. There are many excellent survey articles and books on this topic and generalization of such reconstruction algorithms have been described in the literature (for example, see  \cite{book:Helgason, smith1977,book:tomography01, Quinto-etal} and references therein). Using standard notation  (for example, see \cite{book:Helgason,Quinto06}), the Radon transform $Rf$ of an integrable function $ f $ in $\R^2$ is defined by
\begin{equation}
Rf(\omega,p)=\int_{\langle x,\omega\rangle=p} f(x)dm(x),\label{defRadon1}
\end{equation}
where $\omega=(\omega_1,\omega_2)$ is a unit vector, $p\in \R$, and $dm$ is the arc length measure on the line $\langle x,\omega\rangle =p$ with the usual inner product $\langle~,~\rangle$.
Clearly, the Radon transform can be represented as an integral transform with respect to a measure $ \mu $, which is singular with respect of the Lebesgue measure in $ \R^2 $, as
\begin{equation}
Rf(\omega,p)=\int_{\R^2}  f(x) d \mu = \int_{\R^2} f(x) \delta_{\{\langle x,\omega\rangle =p\}} dx.\label{defRadon2}
\end{equation}
The measure $ \mu $ restricts the Lebesgue measure to lines $E $, parameterized by $ \omega $ and $ p $ in $ \R^2 $ and $\delta_E$ is the Dirac functional on the set $ E $.
%Also, Radon transform is the underlying fundamental concepts for a wide range of other disciplines, including radar imaging, geophysical imaging.

Radon transforms and their inversions are intimately connected to Fourier theory and Riesz potentials and few results deviate from this standard treatment. The standard reconstruction methods are the filtered backprojection method (FBP) and algebraic reconstruction techniques (ART) (for example, see \cite{book:tomography01} for a thorough treatment of these algorithms). Noteworthy among  exceptions to FBP and ART are the works of Milanfar and collaborators \cite{book:Milanfar01, Milanfar02} who use moment-based methods to estimate images from their Radon transform data. \footnote{We became aware of Professor Milanfar's work on moment-based methods for reconstruction of images from its Radon transform data a few weeks after this paper was submitted to this journal. We are happy to have learned about this body of work directly from him and to be able to add references to their work in this paper.} While there are similarities between their results and ours, our main results allow us to reconstruct images from the moments of the modified Radon data. The derived simple relationships between the moments of the modified Radon transform, the moments of the Radon transform and the moments of the objective function appear to be new, and the final inversion algorithm uses these relationships to recover the original function.

If the acquired line integrals are noisy due to statistical fluctuations of photon detection, imperfections of physical system,  or violations in the pencil beam assumption, the classical algorithms are not easy to alter to cancel or reduce the aberrations due to noise. Various authors consider reconstruction algorithms in the presence of noise (for example, see \cite{7123584,Starck:2002:CTI:2319009.2320282}), but our approach is significantly different from those results.

To give a relevant and simple illustration, tomography applications typically collect data in the form of $ Rf $ (instead of $f$) or, in fact $ Rf + \eta $, where $ \eta $ is spurious noise. One viable approach to reduce the effect of the noise is to introduce some mollification to modify the Radon transform after which a sequence of inversions are applied to recover the original density function $f$.
For example, in the context of tomography,  A.K. Louis and P. Maass \cite{Louis01} use projection methods to map their operators into a finite dimensional space determined by the data $ g_N $ and solve the equation $ R f = g_N $.  They assume no knowledge of the inverse of this transform, and using a smoothing operator $ E_\gamma $ (i.e. a convolution operator with a mollifier), with $ \gamma $ being their regularizing parameter, approximate the smoothed (mollified) density function $ f_\gamma = E_\gamma f $.  This is achieved by using $ e_\gamma $ as a suitable mollified basis for the subspace of their Hilbert space,  and approximating $ f_\gamma $ by an element $ v $ in the range of $ R^\ast $ such that $ e_\gamma (x, \cdot) = R^\ast v (x) $. In later works, A. K. Louis extends these ideas to linear operator equations of the first kind \cite{Louis02} and to some nonlinear problems \cite{Louis03}.

In a series of papers Emmanuel Cand\`{e}s and David Donoho (see \cite{candesdonoho00}, and references therein) develop the curvelet transform for reconstruction of images from noisy Radon transforms. These beautiful methods are aimed at detecting edges at certain locations and orientations in the Radon domain and relate these edges to the location and directions of the edges in the original domain.  Their approach makes fundamental use of the fact that curvelets provide an optimal representation of the objects to be identified in the image, thus naturally provide a sparsity in the tight frames used for representing these images and ignores the noise.  While these are sophisticated methods for edge detection,  our goal here has been different and intended to reconstruct the entire image and not just the edges. Furthermore,  to use specific tight frames for general practice is akin to developing an optimal family of mollifiers that isolate particular features of the density function and ignore the noise (e.g. fingerprinting methods).  In designing such mollifiers, the works of Cand\`{e}s and Donoho will be of great interest.

In this paper, we consider the problem of recovering a bivariate moment determinate function from its noisy Radon transform using moments. To reduce the effect of the noise, we introduce a modified Radon transform using mollifiers, and establish an inversion theorem. We derive an explicit relationship between the moments of the modified Radon transform and those of the original function and, using an approximation argument, show that the moment approximations converge uniformly to the original density function. In particular, we recover the original function from the moments of its smoothed Radon transform.
Our strategy differs significantly from the treatments in \cite{Louis01, book:Milanfar01, Milanfar02} since it is based on the moments associated with the modified Radon transform. In addition, the $ L^1 $-methods used in our approach do not lend themselves to Hilbert space projections, and no orthogonality is assumed. As it will be shown the moment integrals are computed directly (as Hamburger moments) in Sections~\ref{sec:MRT} and ~\ref{sec:RMFMRT}, and the approximation in uniform norm to the original density function is derived analytically in Section ~\ref{sec:RDFMND}.

The modified Radon transform is based on convolution of the Radon transform with a symmetric mollifier $\varphi$. The choice of optimal mollifiers for particular applications is an interesting problem that will not be treated in this paper. A simple observation is that the mollifier function may be chosen such that $ (Rf + \eta) \ast \varphi \approx Rf \ast \varphi $.  That is, $ \varphi $ may be chosen such that the noise is washed out or significantly reduced by the mollifier. For example, if the Fourier transforms of $ Rf $ and $ \eta $ have disjoint supports, and $ \varphi $ is chosen such that the support of its Fourier transform is in the support of the Fourier transform of $ Rf $ (or in fact is the same as that support), then $ (Rf + \eta) \ast \varphi = Rf \ast \varphi $. While we cannot in general expect an exact partition of these supports, in many practical applications it is sufficient to reduce the effect of noise on the reconstruction. Since $ \varphi $ can be chosen  to have arbitrarily fast decay outside the support of the Fourier transform of $ Rf $, the portion of the noise spectrally outside the support of the transform of $ Rf $ can be significantly reduced or eliminated.  The question of designing $ \varphi $ in such a way to optimize the recovery of specific features of $ f $ from its noisy Radon transform is an interesting and deep problem. The curvelet transform of Cand\`{e}s and Donoho can be viewed as a special case of this optimization problem, with the mollifiers chosen in such a way to isolate edges in images. %Kapur and Kesavan \cite{Kapur-Kesavan} provide extensive description of use the maximum entropy principle in connection with these optimization methods. 

The remainder of this paper is organized as follows. We begin with collecting pertaining results on Radon transform in Section~\ref{sec:SR}. In Section~\ref{sec:MRT}, we define the modified Radon transform, prove several properties of this transform and derive an inversion formula. Section~\ref{sec:RMFMRT} establishes an explicit relationship between the moments of Radon transform and the moments of the modified Radon transform. In Section~\ref{sec:RDFMND}, we show how to recover $ f $ from the moments of its modified Radon transform and derive theoretical results for the rates of convergence of moment approximations to $ f $. In Section~\ref{sec:ANE}, we present a flow chart that describes the algorithm, and  give a numerical example to demonstrate the convergence of the reconstruction algorithm based on the moments of the mollified Radon transform.  The convergence rates are shown to agree with the estimates derived in Section~\ref{sec:RDFMND}.  Finally, in Section~\ref{sec:concl} we provide a few concluding remarks and summarize the paper.

%
%%%%%%%%%%%%%%%%%%%%%%%%%%%%%%%%%%%%%%%%%%%%%
\section{Standard Results on Radon Transform}
\label{sec:SR}
%%%%%%%%%%%%%%%%%%%%%%%%%%%%%%%%%%%%%%%%%%%%%
For the sake of completeness, this section lists without proof a few well known results about the Radon transform. Readers are referred to references for deeper treatment of these results and their proofs.
A function $f$ is said to be in the \emph{Schwartz space} $\mathcal{S}(\R^2)$ if $f\in C^\infty(\R^2)$ and for each integer $ m \geq 0 $ and each polynomial $P$ of degree $ m $
\begin{equation*}
\underset{x}{\textup{sup}} \big||x|^m P(\partial_{x_1}, \partial_{x_2})f(x)\big|<\infty,
\end{equation*}
where $|x|$ is the Euclidean norm of $x=(x_1,x_2)$. A function $g(\theta,p)$ is said to be in the Schwartz space $\mathcal{S}([0, 2\pi] \times \R)$ if $g(\theta,p)$ can be extended to a smooth and $2\pi$-periodic function in $\theta$, and $g(\cdot,p)\in \mathcal{S}(\R)$ uniformly in $\theta$.  As usual, $ C_c^\infty (\R^2) $ is used to denote $ C^\infty (\R^2) $ functions with compact support.

\begin{lem} \label{lem:range}
For each $f\in \mathcal{S}(\R^2)$, the Radon transform $Rf$ satisfies the following condition: For $k\in \N_0$ the integral
\begin{equation*}
\int_{\R} Rf(\omega,p)p^k~dp
\end{equation*}
is a $k^{th}$ degree homogeneous polynomial in $\omega_1,\omega_2$ \textup{(see \cite{book:Helgason}, Lemma 2.3)}.
\end{lem}

We denote the unit vector in direction $\theta$ as $\omega=\omega(\theta):=(\omega_1,\omega_2)$ with $\omega_1=\cos{\theta}$ and  $\omega_2=\sin{\theta}$. Thus, the Radon transform of $f\in L^1(\R^2)$ can be expressed as a function of $(\theta,p)$:
\begin{equation}
Rf(\theta,p)=\int_{\langle x,\omega(\theta)\rangle=p} f(x)dx.\label{defRadon3}
\end{equation}

Note that since the pairs $(\omega,p)$ and $(-\omega,-p)$ give the same line, $R$ satisfies the \emph{evenness} condition: $Rf(\theta,p)=Rf(\theta+\pi,-p)$.

% Radon transform R is continuous
\begin{thm}\label{radon:cont}
The Radon transform $R$ is a bounded linear operator from $L^1(\R^2)$ to $L^1([0,2\pi]\times \R)$
with norm $\|R\|\leq 2\pi$, i.e., $\|R f\|_{L^1([0,2\pi]\times \R)} \leq 2\pi \|f\|_{L^1(\R^2)}$.
\begin{proof}
See \cite{Quinto06}, for example.
\end{proof}
\end{thm}

Along with the transform $Rf$, we define the dual Radon transform of $g\in L^1([0,2\pi] \times \R)$ as
\begin{equation}
R^*g(x)=\int_{0}^{2\pi} g(\theta, \langle x, \omega \rangle) d\theta,
\end{equation}
which is the integral of $g$ over all lines that go through $x$.

Using $F_1$ and $F_2$ for the 1-D and 2-D Fourier transforms of integrable functions with integrable transforms, recall that
\begin{align*}
F_1f(s) &= \frac{1}{\sqrt{2\pi}} \int_{-\infty}^{\infty} f(t) e^{-ist} dt,\\
F_1^{-1} f(t) &= \frac{1}{\sqrt{2\pi}} \int_{-\infty}^{\infty} f(s) e^{ist} ds,\\
F_2f(\xi) &= \frac{1}{2\pi} \int_{\R^2} f(x) e^{-i\langle x,\xi\rangle } dx,\\
F_2^{-1}f(x) &= \frac{1}{2\pi} \int_{\R^2} f(\xi) e^{i\langle x ,\xi \rangle } d\xi.
\end{align*}

% Projection-slice Theorem
\begin{thm}[Projection-Slice Theorem \cite{book:Helgason},\cite{Quinto06}]\label{thm:projection-slice}
Let $f\in L^1(\R^2)$.
Then,
\begin{equation*}
F_2 f(s\omega)=\frac{1}{\sqrt{2\pi}} F_1 (Rf(\theta,\cdot))(s).
\end{equation*}
\end{thm}
This theorem shows that $R$ is injective on $L^1(\R^2)$.
The Fourier inversion formula combined with the Projection-Slice Theorem provides an inversion formula for the Radon transform in $ \R^2 $.

% Riesz Potential
Denote the Riesz potential $I^{-1}$, for $g\in L^1([0,2\pi]\times \R)$, as the operator with Fourier multiplier
$|s|$ (see \cite{book:Stein}):
\begin{equation}
I^{-1}g \ (\theta, t) =F_1^{-1}(|s|(F_1g(\theta,t)(s)) )(t).
\end{equation}

% Inversion formula for R
\begin{thm}[Inversion formula for $R f$ \cite{book:Helgason,Quinto06}]\label{thm:inversionf}
Let $f\in C_c^\infty(\R^2)$. Then
\begin{equation*}
f(x)=\frac{1}{4\pi}R^* (I^{-1} Rf)(x).
\end{equation*}
\end{thm}
Note that this theorem is true on a larger domain than $C_c^\infty(\R^2)$. However, $I^{-1}Rf$ may be a distribution rather than a function.
%

% Question
One may ask, when is a given function $g$ the Radon transform of a function $f$? In other words, for a given function $g$, does there exist $f$ such that $g=Rf$? The following theorem, which partly motivates the naturality of use of moments in this context, is the fundamental result on this question. This result is called the \emph{Schwartz or Range theorem} for the Radon transform and its proof is established in \cite{book:Helgason}, Theorem 2.4.
%
% Range Theorem
\begin{thm}
Let $g\in \sw(\domain)$ be even. Then, there exists $f\in \sw(\R^2)$ such that $g=Rf$
if and only if for each $k\in \N_0$, the $k^{th}$-moment
\begin{equation} \label{homo}
\int_{-\infty}^{\infty} g(\theta,p)p^kdp
\end{equation}
is a homogeneous polynomial of degree $k$ in  $\omega_1$ and $\omega_2$.
\end{thm}

%
%%%%%%%%%%%%%%%%%%%%%%%%%%%%%%%%%%%%%%%%%%%%%%%%%%
\section{Modified Radon Transform} \label{sec:MRT}
%%%%%%%%%%%%%%%%%%%%%%%%%%%%%%%%%%%%%%%%%%%%%%%%%%

% Introduction to mollifier
\begin{defn}\label{defn:molifier}
If $\varphi$ is a smooth, nonnegative, and compactly supported integrable function on $\R$, such that
\begin{enumerate}
\item ${\displaystyle \int_{-\infty}^{\infty} \varphi(t)dt=1}$,
\item Setting $\varphi _{h}(t) :=  \frac{1}{h} \varphi (\frac{t}{h})$, and defining the convolution of $ f $ and $ \varphi_h $ by
$$ 
f \ast \varphi_h (x) = \int_{\mathbb R} f(y) \varphi_h (x-y) dy = \int_{\mathbb R} f(x-y) \varphi_h (y) dy , 
$$
then $\|f \ast \varphi_h - f\|_{L^\infty(\mathbb{R})} \rightarrow 0$ as $h \rightarrow 0$.
\end{enumerate}
Such a $ \varphi $ is called a \textit{(positive) mollifier}. Furthermore, if
$\varphi(t) = g (|t|)$ for some infinitely differentiable function
$g :{\R}\rightarrow \R$, then $ \varphi $ is called a \textit{symmetric mollifier}.
\end{defn}

% Example of a molifier
For example, if $\varphi:\R\rightarrow \R$ is { defined as follows}
%the Gaussian function
\begin{equation*}
\varphi(x) = \left\{
\begin{array}{ll}
Ae^{-\frac{1}{1-|x|^2}} & \text{if } |x| < 1\\
0 & \text{if } |x| \geq 1,
\end{array} \right.
\end{equation*}
then $\varphi$ is a positive symmetric mollifier. Here $A$ is a constant such that (i) in Definition~\ref{defn:molifier} is satisfied.

% Definition of modified Radon transform
Let $\Omega=\{\varphi : \varphi$ is a positive symmetric mollifier such that $F_1(\varphi)(s)>0$ for all $s\in\R \}$. Clearly, the Gaussian function is in $\Omega$. 
Since $ \varphi $ is smooth and compactly supported, it has finite moments for all orders.

\begin{defn}\label{defn:MRT}
Let $\varphi\in \Omega$ and $f\in L^1(\R^2)$.
The \emph{modified Radon transform} in $\R^2$ is defined by
\begin{equation}
\widehat{R}_\varphi f(\theta,p) = \int_{\R^2} \chi(x;\theta,p) f(x)dx,
\end{equation}
where
\begin{equation*}
\chi(x;\theta,p) = (\delta \ast \varphi) (\langle x,\omega \rangle-p),~ (\delta \textup{ is a delta function}).
\end{equation*}
\end{defn}

Next we show that the modified Radon transform, which is defined as the smoothed Radon transform of the original density is equal to the Radon transform of the smoothed density.  

% Lemma
\begin{prop}\label{modified:form2}
If $ \varphi \in \Omega $ and $ f \in L^1 ({\mathbb R}^2) $, then
\begin{equation}
\widehat{R}_\varphi f(\theta,p) = \int_{\R} Rf(\theta,p+\tau) \varphi (\tau) d\tau = \Big(Rf(\theta,\cdot) \ast \varphi \Big)(p).
\end{equation}
\begin{proof}
\begin{align*}
\widehat{R}_\varphi f(\theta,p) & = \int_{{\mathbb R}^2} (\delta \ast \varphi)(\langle x,\omega\rangle-p) f(x)dx\\
& = \int_{\mathbb R} \Bigg( \int_{{\mathbb R}^2} \delta(\langle x,\omega \rangle-p-\tau) f(x)dx \Bigg) \varphi (\tau) d\tau \\
& = \int_{{\mathbb R}} \Bigg( \int_{\langle x,\omega \rangle=p+\tau} f(x)dx \Bigg) \varphi (\tau) d\tau \\
%\end{align*}
%\end{proof}
%\end{lem}
% Proposition
%\begin{prop}\label{modified:form2}
%$$\begin{equation*}
%\widehat{R}_\varphi f(\theta,p) = & = %\Big(Rf(\theta,\cdot) \ast \varphi \Big)(p).
%\end{equation*}
%\begin{proof}
%By Lemma \ref{modified:form1}
%\begin{align*}
%\widehat{R}_\varphi f(\theta,p)
& = \int_{\mathbb R} Rf(\theta,p+\tau) \varphi (\tau) d\tau \\
& = \int_{\mathbb R} Rf(\theta,p-\tau) \varphi (-\tau) d\tau.
\end{align*}
Since $\varphi$ is a symmetric mollifier,
\begin{equation*}
\widehat{R}_\varphi f(\theta,p)
= \int_{\R} Rf(\theta,p-\tau) \varphi (\tau) d\tau.
\end{equation*}
\end{proof}
\end{prop}

% Theorem (modified Radon transform is a bounded operator)
\begin{thm}\label{boundedL1}
The modified Radon transform $\widehat{R}_\varphi$ is a bounded linear operator from $L^1(\R^2)$ to $L^1([0,2\pi]\times \R)$ with norm $\|\widehat{R}_\varphi\|\leq 2\pi$,i.e.,
$\|\widehat{R}_\varphi f\|_{L^1([0,2\pi]\times \R)} \leq 2\pi \|f\|_{L^1(\R^2)}$.
\begin{proof}
Using Proposition \ref{modified:form2}, one finds that
\begin{align*}
\|\widehat{R}_\varphi f\|_{L^1([0,2\pi] \times \R)}
&= \int_{\theta=0}^{2\pi}  \int_{p=-\infty}^{\infty}  \Big|\widehat{R}_\varphi f(\theta,p)\Big|dp d\theta\\
&= \int_{\theta=0}^{2\pi}  \int_{p=-\infty}^{\infty}
\Bigg|\int_{\tau=-\infty}^{\infty} Rf(\theta, p+\tau) \varphi (\tau) d\tau \Bigg| dp d\theta\\
&\leq \int_{\theta=0}^{2\pi}  \int_{p=-\infty}^{\infty}
\int_{\tau=-\infty}^{\infty} \Big|Rf(\theta, p+\tau)\Big| \varphi (\tau) d\tau  dp d\theta\\
&= \int_{\tau=-\infty}^{\infty} \varphi(\tau) \Bigg( \int_{\theta=0}^{2\pi}  \int_{p=-\infty}^{\infty}
 \Big|Rf(\theta, p+\tau)\Big| dp d\theta \Bigg)  d\tau.
\end{align*}
By Theorem \ref{radon:cont} and Definition \ref{defn:molifier}, it follows that
\begin{align*}
\|\widehat{R}_\varphi f\|_{L^1([0,2\pi] \times \R)}
&\leq \int_{\tau=-\infty}^{\infty} \varphi(\tau) \Big( 2\pi \|f\|_{L^1(\R^2)} \Big)  d\tau\\
&= 2\pi \|f\|_{L^1(\R^2)}.
\end{align*}
\end{proof}
\end{thm}

Denote the modified Riesz potential $\widehat{I}^{-1} $, for $g\in L^1([0,2\pi]\times \R)$, as the operator with Fourier multiplier $|s|$ and symmetric mollifier $\varphi \in \Omega $:
\begin{equation*}
\widehat{I}^{-1}g (\theta) =\frac{1}{\sqrt{2\pi}} F_1^{-1}\bigg(|s|\Big(\frac{F_1(g(\theta,\cdot))(s)}{F_1(\varphi) (s)}\Big)\bigg).
\end{equation*}

% Theorem (Inversion formula for modified Radon transform)
The proof of Theorem~\ref{thm:inversionf} combined with the convolution theorem
provides an inversion formula for $f$ from $\widehat{R}_\varphi$.
\begin{thm}[Inversion formula for $\widehat{R}_\varphi f$]\label{inv-modified}
Let $f\in C_c^{\infty}(\R^2)$. Then
\begin{equation*}
f(x)=\frac{1}{4\pi}R^* (\widehat{I}^{-1} \widehat{R}_\varphi f)(x).
\end{equation*}
\begin{proof}
By Proposition \ref{modified:form2} and the convolution theorem,
\begin{equation}\label{conv}
F_1(\widehat{R}_\varphi f(\theta,\cdot))(s)
= \sqrt{2 \pi} F_1(Rf(\theta,\cdot))(s) F_1 (\varphi)(s).
\end{equation}
Since $F_1 (\varphi)(s)\neq 0$, by Theorem \ref{thm:projection-slice},
it follows that
\begin{equation}
F_2f(s\omega) = \frac{1}{2\pi} \frac{F_1(\widehat{R}_\varphi f(\omega,\cdot))(s)}{F_1 (\varphi)(s)}. \label{eq:mod-proj-slice}
\end{equation}
Applying the Fourier inversion formula and Theorem \ref{thm:inversionf}, it follows that
\begin{align*}
f(x) &= \frac{1}{2\pi} \int_{\R^2} F_2f(\xi) e^{i<x,\xi>} d\xi\\
&= \frac{1}{2\pi}\frac{1}{2} \int_{\theta=0}^{2\pi}
\int_{s=-\infty}^{\infty} F_2f(s\omega) e^{i<x,s\omega>} |s| ds d\theta\\
&= \frac{1}{4\pi} \int_{\theta=0}^{2\pi}
\widehat{I}^{-1}\widehat{R}_\varphi f(\theta,\langle x, \omega \rangle)d\theta\\
&=\frac{1}{4\pi} R^* (\widehat{I}^{-1} \widehat{R}_\varphi f)(x).
\end{align*}
\end{proof}
\end{thm}

%
%%%%%%%%%%%%%%%%%%%%%%%%%%%%%%%%%%%%%%%%%%%%%%%%%%%%%%%%%%%%%%
\section{Recovering Moments from the Modified Radon Transform}
\label{sec:RMFMRT}
%%%%%%%%%%%%%%%%%%%%%%%%%%%%%%%%%%%%%%%%%%%%%%%%%%%%%%%%%%%%%%

%Suppose that a cumulative distribution function (cdf) $F$ is absolutely continuous with respect to the Lebesgue measure with support $[0,1]^2$. Denote the corresponding density function by $f$. 

Let $\omega=(\cos{\theta},\sin{\theta})$ denote a unit direction vector and $x=(x_1, x_2)$ a vector in $\R^2$.
Suppose that $f:\R^2 \rightarrow \R$ is in the Schwartz space. Then by Lemma~\ref{lem:range}, the definition of Radon transform and Fubini's theorem, we have
\begin{equation}\label{range}
\int_{-\infty}^{\infty} Rf(\theta,p)p^k~dp = \int_{\R^2} f(x)\langle \omega, x \rangle^k dx \quad \text{for each }k\in \N_0.
\end{equation}
Appropriate expansion of the right hand side of \eqref{range} using definition of moment gives
%Assume that the support of $f$ is contained within $I^2$ with $I=[0,1]$.
%Then the identity \eqref{range} can be expressed as
\begin{equation}\label{system1}
b^{(k)}(\theta) := \int_{-\infty }^{\infty } Rf(\theta,p) p^k~dp =\sum_{j=0}^k
C(k,j)
\big(\cos^{j}{\theta}\sin^{k-j}{\theta}\big)\gamma_{j,k-j} (f),
\end{equation}
where 
\begin{equation*}
C(k,j)=\dfrac{k!}{j!(k-j)!}
%\begin{pmatrix}
%k\\
%j\\
%\end{pmatrix}
\quad \text{and} \quad
\gamma_{\alpha_1,\alpha_2} (f)=\int_{\R^2} x_1^{\alpha_1} x_2^{\alpha_2} f(x)dx,~ \alpha_1, \alpha_2 \in \N_0.
\end{equation*}

Let $0< \theta_0 < \theta_{1} < \cdots < \theta_{k} < \pi$ be distinct angles.
A sampling of \eqref{system1} on these angles yields a linear algebraic form of dimension $k+1$ written as
\begin{equation}\label{axb2}
\mathbf{A}^{(k)}\mathbf{x}^{(k)}=\mathbf{b}^{(k)},
\end{equation}
where 
\small
\begin{equation*}
\mathbf{A}^{(k)}=
\begin{pmatrix}
C(k,0)\cos^0{\theta_0}\sin^k{\theta_0} & C(k,1)\cos^1{\theta_0}\sin^{k-1}{\theta_0} & \cdots & C(k,k)\cos^k{\theta_0}\sin^0{\theta_0}\\
C(k,0)\cos^0{\theta_1}\sin^k{\theta_1} & C(k,1)\cos^1{\theta_1}\sin^{k-1}{\theta_1} & \cdots & C(k,k)\cos^k{\theta_1}\sin^0{\theta_1}\\
C(k,0)\cos^0{\theta_2}\sin^k{\theta_2} & C(k,1)\cos^1{\theta_2}\sin^{k-1}{\theta_2} & \cdots & C(k,k)\cos^k{\theta_2}\sin^0{\theta_2}\\
\vdots & \vdots & \ddots & \vdots \\
C(k,0)\cos^0{\theta_k}\sin^k{\theta_k} & C(k,1)\cos^1{\theta_k}\sin^{k-1}{\theta_k} & \cdots & C(k,k)\cos^k{\theta_k}\sin^0{\theta_k}
\end{pmatrix},
\end{equation*}
\normalsize
\begin{equation*}
\mathbf{x}^{(k)}=
\begin{pmatrix}
\gamma_{0,k}(f)\\
\gamma_{1,k-1}(f)\\
\gamma_{2,k-2}(f)\\
\vdots\\
\gamma_{k,0}(f)\\
\end{pmatrix},\quad 
\mathbf{b}^{(k)}=
\begin{pmatrix}
b^{(k)}(\theta_0)\\
b^{(k)}(\theta_1)\\
b^{(k)}(\theta_2)\\
\vdots\\
b^{(k)}(\theta_k)\\
\end{pmatrix}.
\end{equation*}
The determinant of the matrix $\mathbf{A}^{(k)}$ can be expressed as:
\begin{equation*}
\det(\mathbf{A}^{(k)})=\det(V^{(k)}) \prod_{j=1}^k \big( C(k,j) \sin^k{\theta_j} \big) ,
\end{equation*}
where $V^{(k)}=[\cot^{j-1}{\theta_i}]_{1\leq i,j \leq k+1}$ is a Vandermonde matrix.
Using the Vandermonde determinant formula, it is easy to show $\det(\mathbf{A}^{(k)})$ is positive, implying
the system (\ref{axb2}) has a unique solution $\mathbf{x}^{(k)}$.
Note that the matrix $\mathbf{A}^{(k)}$ is positive definite since its leading principal minors are all positive.

% Theorem (Range Theorem for modified Radon transform)
Using the above, we may establish an algebraic relation between  the moments of the modified Radon transform and the moments of the Radon transform.
\begin{thm}\label{modified:range1}
Let $ f \in \mathcal{S}(\R^2)$. Then for each $k\in \N_0$,
\begin{equation}\label{modified:range2}
\hat{b}^{(k)}(\theta):=\int_{-\infty}^{\infty} \widehat{R}_\varphi f(\theta,p)p^kdp
= \sum_{j=0}^k C(k,j) c_{j}
b^{(k-j)}(\theta)
\end{equation}
where $\displaystyle{c_{j}=(-1)^j\gamma_{j} (\varphi)}$ and $\displaystyle \gamma_{j}(\varphi)=\int_{-\infty}^{\infty}\tau^j  \varphi(\tau)d\tau$
for each $j\in \N_0$. 
That is, the value $\hat{b}^{(k)}(\theta)$ is a linear combination of $b^{(0)}(\theta)$, $b^{(1)}(\theta)$, $\ldots$, $b^{(k)}(\theta)$ for any $\theta$.
\begin{proof}
By Proposition \ref{modified:form2}, it follows that
\begin{align*} % homogeneous polynomial
\int_{-\infty}^{\infty} \widehat{R}_\varphi  f(\theta,p)p^k~dp
&= \int_{-\infty}^{\infty}
\Bigg(\int_{-\infty}^\infty Rf(\theta,p+\tau) \varphi (\tau) d\tau \Bigg) p^k dp\\
&= \int_{-\infty}^{\infty}  \int_{-\infty}^\infty
\Bigg(\int_{\langle x,\omega\rangle =p+\tau} f(x)dx \Bigg) \varphi (\tau) d\tau~ p^k dp\\
&= \int_{-\infty}^\infty \varphi (\tau) \Bigg[ \int_{-\infty}^{\infty} p^k
\Bigg(\int_{\langle x,\omega\rangle =p+\tau} f(x)dx \Bigg) dp \Bigg] d\tau\\
&= \int_{-\infty}^\infty \varphi (\tau)
\Bigg(\int_{\R^2} (\langle x,\omega\rangle -\tau)^k f(x)dx \Bigg) d\tau.\\
\end{align*}
Using the following polynomial expansion
\begin{align*}
(\langle x,\omega\rangle-\tau)^k = \sum_{|\alpha|=k} \frac{k!}{\alpha!}
(x_1\cos{\theta})^{\alpha_1}(x_2\sin{\theta})^{\alpha_2}(-\tau)^{\alpha_3},
\end{align*}
where $\alpha=(\alpha_1,\alpha_2,\alpha_3)$ is a multi-index with $|\alpha|:= \alpha_1 + \alpha_2 + \alpha_3$ and $\alpha!:=\alpha_1!\alpha_2!\alpha_3!$, it follows that
\begin{align*}
\int_{-\infty}^{\infty} {\widehat R}_\varphi f(\theta,p)p^k dp\\
&\hspace*{-2.5cm}= \sum_{|\alpha|=k}\frac{k!}{\alpha!}c_{\alpha_3}
\gamma_{\alpha_1,\alpha_2}(f) \cos^{\alpha_1}{\theta} \sin^{\alpha_2}{\theta}\\
&\hspace*{-2.5cm}= \sum_{\alpha_3=0}^k C(k,\alpha_3) c_{\alpha_3}
\sum_{j=0}^{k-\alpha_3}
C(k-\alpha_3,j)
\gamma_{j,k-\alpha_3-j}(f) \cos^{j}{\theta} \sin^{k-\alpha_3-j}{\theta}.
\end{align*}
\end{proof}
\end{thm}
Let $0< \theta_1 < \theta_{2} < \cdots < \theta_{k} < \pi$ be distinct angles. A sampling of \eqref{modified:range2} gives
\begin{equation}\label{bbc}
%\mathbf{\widehat{B}}^{(k)} = \mathbf{B}^{(k)}\mathbf{C}^{(k)},
\mathbf{C}^{(k)} \mathbf{B}^{(k)} = \mathbf{\widehat{B}}^{(k)},
\end{equation}
where
\begin{equation*}
\mathbf{C}^{(k)}=
\begin{pmatrix}
C(k,0)c_0 & 0 & \cdots & 0\\
C(k,1)c_1 &  C(k-1,0)c_0 & \cdots & 0\\
C(k,2)c_2 & C(k-1,1)c_1 & \cdots & 0\\
\vdots & \vdots & \ddots & \vdots \\
C(k,k)c_k & C(k-1,k-1)c_{k-1} & \cdots & C(0,0)c_0\\
\end{pmatrix}^\top,
\end{equation*}
\begin{equation*}
\mathbf{B}^{(k)}=
\begin{pmatrix}
b^{(k)}(\theta_0) & b^{(k-1)}(\theta_0) & \cdots & b^{(0)}(\theta_0)\\
b^{(k)}(\theta_1) & b^{(k-1)}(\theta_1) & \cdots & 0\\
b^{(k)}(\theta_2) & b^{(k-1)}(\theta_2) & \cdots & 0\\
\vdots & \vdots & \ddots & \vdots \\
b^{(k)}(\theta_k) & 0 & \cdots & 0
\end{pmatrix}^\top,
\end{equation*}
\begin{equation*}
\mathbf{\widehat{B}}^{(k)}=
\begin{pmatrix}
\hat{b}^{(k)}(\theta_0) & \hat{b}^{(k-1)}(\theta_0) & \cdots & \hat{b}^{(0)}(\theta_0)\\
\hat{b}^{(k)}(\theta_1) & \hat{b}^{(k-1)}(\theta_1) & \cdots & 0\\
\hat{b}^{(k)}(\theta_2) & \hat{b}^{(k-1)}(\theta_2) & \cdots & 0\\
\vdots & \vdots & \ddots & \vdots \\
\hat{b}^{(k)}(\theta_k) & 0 & \cdots & 0
\end{pmatrix}^\top.
\end{equation*}

Since the matrix $\mathbf{C}^{(k)}$ is upper triangular, the determinant of the matrix is $(c_0)^{k+1}>0$.
Then there exists a unique solution $\mathbf{B}^{(k)}$ for the given matrix $\mathbf{\widehat{B}}^{(k)}$.
Thus we have the vectors $\mathbf{b}^{(i)}$ for all $i=0,1,\ldots,k$. 
By the equation \eqref{axb2}
one can find $\mathbf{x}^{(0)}$, $\mathbf{x}^{(1)}$, \ldots, $\mathbf{x}^{(k)}$ such that
$\mathbf{x}^{(i)}=(\mathbf{A}^{(i)})^{-1}\mathbf{b}^{(i)}$ for all $i=0, 1,\ldots, k$. That is,
one has the moments $\{\gamma_{\alpha_1,\alpha_2}(f) \}_{\alpha_1+\alpha_2 \leq k}$.

%
%
%%%%%%%%%%%%%%%%%%%%%%%%%%%%%%%%%%%%%%%%%%%%%%%%%%%%%%%%%%%%%%%%%%
\section{Recovery of Density Function via Moments with Noisy Data}
\label{sec:RDFMND}
%%%%%%%%%%%%%%%%%%%%%%%%%%%%%%%%%%%%%%%%%%%%%%%%%%%%%%%%%%%%%%%%%%
By 
%Lemma \ref{modified:form1} and 
{Proposition \ref{modified:form2}}, since the Radon transform of a smoothed function is the smoothed Radon transform of that  function, the moments of the Radon transform of $ f $ and the modified Radon transform may be related by combining Theorem
\ref{modified:range1} and equation \eqref{system1}.

Recall that the convolution $f\ast \varphi$ of $f$ and $\varphi$ is defined as follows:
\begin{equation*}
(f\ast \varphi) (x)=\int f(x_1 -\tau, x_2 -\tau) \varphi (\tau)\, d \tau \;\;\; {\rm for \;\; each}\;\;\; x=(x_1, x_2).
\end{equation*}

\begin{thm}
\label{synthesis}
If $ f \in \mathcal{S}(\R^2)$ and the density function $ f $ is recovered from the modified Radon transform, the moments must satisfy a necessary linear constraint given by
$$
\displaystyle{
\begin{array}{l}
\displaystyle
\sum_{j=0}^k \sum_{\ell=0}^{j} C(k,j)  C(j, \ell)  \\
\hspace{0.3in} \Big[
\gamma_{k-j} (\varphi) \gamma_{\ell,j-\ell} (f) \cos^\ell \theta \sin^{j-\ell} \theta \\
\hspace{0.3in} - \displaystyle\sum_{n=0}^{k-j}  C(k-j,n) \gamma_{j-\ell, k-j-n} (f) \gamma_{n+\ell} (\varphi) \cos^j \theta \sin^{k-j} \theta\Big] = 0. \\
\end{array}
}
$$

\begin{proof}
By Proposition \ref{modified:form2},  
$ \widehat{R}_\varphi f(\theta,p) = \Big(Rf(\theta,\cdot) \ast \varphi \Big)(p) $. The $ k $-th moment of the modified Radon transform on the left hand side of this equation is given by
\begin{equation}\label{eq:4-mrt2}
\begin{aligned}
\displaystyle &\int_{-\infty}^{\infty} \widehat{R}_\varphi f(\theta,p)p^k dp 
%\\&
= \displaystyle\sum_{j=0}^k C(k,j) \cos^j \theta \sin^{k-j} \theta \gamma_{j,k-j} (f \ast \varphi)\\ 
 &=\displaystyle\sum_{j=0}^k \sum_{l=0}^{j} C(k,j) C(j,\ell)
 \displaystyle\sum_{n=0}^{k-j} C(k-j, n) \gamma_{j-\ell, k-j-n} (f) \gamma_{n+\ell} (\varphi) \cos^j \theta \sin^{k-j}\theta,
\end{aligned}
\end{equation}

%\begin{align}
%\label{eq:4-mrt2}
%\displaystyle{
%\begin{array}{l}
%\int_{-\infty}^{\infty} \widehat{R}_\varphi f(\theta,p)p^kdp\\ \notag 
%\hspace{0.3in} = \sum_{j=0}^k C(k,j) \cos^j \theta \sin^{k-j} \theta \gamma_{j,k-j} (f \ast \varphi)\notag \\
%\hspace{0.3in} = \sum_{j=0}^k \sum_{l=0}^{j} C(k,j) C(j,\ell) \notag \\
%\hspace{0.3in}= \sum_{n=0}^{k-j} C(k-j, n) \gamma_{j-\ell, k-j-n} (f) \gamma_{n+\ell} (\varphi) \cos^j \theta \sin^{k-j}\theta ,
%\end{array}
%}
%\end{align}
% \nonumber\\
while the $k$-th moment of the right hand side is given by
$$
\begin{array}{l}
\displaystyle\int_{-\infty}^{\infty} \Big(Rf(\theta,\cdot) \ast \varphi \Big)(p) p^kdp  \\
\hspace{0.3in}= \displaystyle\int_{-\infty}^{\infty} \left[ \int_{-\infty}^{\infty} Rf(\theta, p - \tau) \varphi(\tau) d \tau \right] p^k dp \\
\hspace{0.3in}= \displaystyle\int_{-\infty}^{\infty} \varphi(\tau) \left[ \int_{-\infty}^\infty Rf(\theta, u) (u + \tau)^k du \right] d\tau \\
\hspace{0.3in} = \displaystyle \sum_{j=0}^k C(k,j) \left[ \int_{-\infty}^{\infty} \varphi (\tau) \tau^{k-j} d \tau \right] \left[\int_{-\infty}^{\infty} Rf(\theta, u) u^j du \right] \\
\hspace{0.3in} = \displaystyle\sum_{j=0}^k \sum_{\ell=0}^{j} 
%(-1)^{j} 
C(k,j) C(j, \ell) \gamma_{k-j} (\varphi) 
\gamma_{\ell,j-\ell}(f)  \cos^{\ell}{\theta} \sin^{j-\ell}{\theta}.
\end{array}
$$

Equating these expressions and rearranging, the expression in the statement of the theorem follows.
\end{proof}
\end{thm}

It is noteworthy that Theorem \ref{synthesis} relates the moments of the modified Radon transform to the moments of the density function $ f $. Hence the inversion theorem (Theorem \ref{inv-modified}) is implicitly the inversion of this system of equations. 

\medskip

 Now suppose that $f$ is a moment determinate measurable function with compact support $[0,1]^2$ such that $ \int_{[0,1]^2} f = 1 $. 
 Note that a function is moment determinate if it is uniquely determined from its moments. 
 For an authoritative discussion of moment determinate functions see \cite{book:Akhiezer}. Let $\{\gamma_{\alpha_1,\alpha_2}(f) \}_{\alpha_1\leq m,\alpha_2 \leq n}$ be a given sequence of moments of $f$ up to {order} $m+n$. In \cite{ML13}, R. Mnatsakanov and S. Li construct the approximation of $f$, denoted by $ app(f) $,  using the moments of $ f $ of order up to 
$ m + n$  as 
\begin{equation}\label{M-Li}
C_{m,n}(x) \sum_{\alpha_1=0}^{m-[mx_1]} \sum_{\alpha_2=0}^{n-[nx_2]}
\frac{(-1)^{\alpha_1+\alpha_2} \gamma_{\alpha_1 + [mx_1],\alpha_2 + [nx_2]}(f) }{\alpha_1!\alpha_2!(m-[mx_1]-\alpha_1)!(n-[nx_2]-\alpha_2)!}
\end{equation}
where 
$$ 
C_{m,n} (x) = \frac{\Gamma(m+2)\Gamma(n+2)}{\Gamma([m x_1] +1) \Gamma([n x_2] +1)} 
$$ 
for all $\alpha_1, \alpha_2\in\N_0$ with $\alpha_1\leq m$ and $\alpha_2 \leq n$. {Here by $[a]$ we denote the integer part of $a$.}

In the following we show that the approximation of $ f $ given by equation (\ref{M-Li}) converges to $ f $ uniformly, and provide a rate for this convergence.  To do so,  let $ f_{10} $, $ f_{01} $, $f_{20} $, $f_{11}$, $ f_{02}$  denote the partial derivatives of $ f $ up to the second order.  We have

\begin{thm}
%\cite{ML13}
\label{approx1}
Let $ f \in C^2( [0,1]^2) $.  Then
%$f\in C^{2} ([0,1]^2)$. Then 
$app( f ) $ converges uniformly to $ f $  on $ [0,1]^2$ as $ m, n \rightarrow \infty $,  and
\begin{align}\label{eq:4-1}
||app(f) - f||
& \leq  \frac{ 2\,  || f^{}_{10}|| + \frac{1}{2}\, ||f^{}_{20}|| }{m+2} +  \frac{2\,   ||f^{}_{01}|| + \frac{1}{2}\, ||f^{}_{02}||}{n+2} \notag \\
&+  \frac{\frac{1}{2}\, ||f^{}_{11}||}{\sqrt{(m +2)\,(n+2})}
+ o\Big (\frac{1}{m}\Big ) + o\Big (\frac{1}{n}\Big )+ o\Big (\frac{1}{\sqrt{m\, n}}\Big ). 
\end{align}
In particular, by choosing $m=n$, we have,  as $ n \rightarrow \infty$:
\begin{align}\label{eq:4-2}
||app(f) - f||
& \leq  \frac{C}{n+2} + o\Big (\frac{1}{n}\Big ),
\end{align}
 %as $ n \rightarrow \infty$.
 %Here $C= 2 \Big (||f^{}_{10}|| +  ||f^{}_{01}||\Big )  + \frac{1}{2}\, \Big ( ||f^{}_{20}|| + || f^{}_{11}|| +  ||f^{}_{02}|| \Big )$.
%\end{thm}
%Assume that the density function $f$ is continuous on $[0,1]^2$. Then $app(f) $ converges uniformly to $ f $ as $m,n \rightarrow \infty$. Furthermore, for any $\delta>0$
%\begin{equation*}\label{approx}||app(f) - f || \leq \Delta(f,\delta) + \frac{4||f||}{\delta^2(\alpha^* +2)}+\frac{2||f||}{\delta^4 (m+2)(n+2)},
%%\end{equation*}
where $||\cdot||$ is the sup-norm and $C= 2 \Big (||f^{}_{10}|| +  ||f^{}_{01}||\Big )  + \frac{1}{2}\, \Big ( ||f^{}_{20}|| + || f^{}_{11}|| +  ||f^{}_{02}|| \Big )$.
% $\alpha^*=\min{\{m,n\}}$, and $\Delta(f,\delta)$ represents the modulus of continuity of $f$. 
%\end{thm}

\begin{proof} Let us plug in the moments of $f$:
\begin{align}\label{eq:4-3}
\gamma_{\alpha_1, \alpha_2} (f) = \int_{0}^{1} \int_{0}^{1} t^{\alpha_1}\, s^{\alpha_2} f(t, s){d}t \,{d}s
\end{align}
into the formula that defines $app(f)$, Eq. (\ref{M-Li}).  After 
%changing the order of summations and  integration,  and 
applying Newton's binomial formula twice, we easily obtain:
\begin{align*}
&\sum_{\alpha_1=0}^{m-[mx_1]} \sum_{\alpha_2=0}^{n-[nx_2]}
\frac{(-1)^{\alpha_1+\alpha_2}  t^{\alpha_1}\, s^{\alpha_2} }{\alpha_1!\alpha_2!(m-[mx_1]-\alpha_1)!(n-[nx_2]-\alpha_2)!}\notag \\
&=\frac{ (1-t)^{m-[mx_1]}}{(m-[mx_1])!} \ \ 
\frac{(1-s)^{n-[nx_2]}  }{(n-[nx_2])!} .
\end{align*}
Hence, using these two steps combined with changing the order of summations and  integration, one obtains:
\begin{align} \label{eq:4-4}
\hspace*{-0.35cm} app(f) (x)- f (x)&=\int_{0}^{1} \int_{0}^{1} f (t, s)   \beta_{m, x_1}(t)  \beta_{n, x_2}(s)  d t\, ds -  f (x) \notag \\
& = \int_{0}^{1} \int_{0}^{1} \Big (f (t, s) - f (x_1, x_2) \Big ) \beta_{m, x_1}(t)  \beta_{n, x_2}(s)  d t\, ds.
\end{align}
Here by $ \beta_{m, u}(\cdot):=\beta(\cdot, [m u] +1, m - [m u]+ 1)$ we denote the beta density function with shape parameters $ [m u]+1$ and $m - [m u]+ 1$. Note that the mean and  variance of $ \beta_{m, u}(\cdot)$ are:
\begin{align}\label{eq:4-5}
\theta_{m, u}=\frac{[m u] +1}{m+2} \;\;\;\;\; {\rm and }\;\;\;  \sigma^2_{m, u}=\frac{([m u] +1)( m - [m u]+ 1)}{(m+2)^2\, (m+3)},
\end{align}
respectively. Also, it is worth mentioning that the sequence  of functions $\{ \beta_{m, x}(t): m \geq 1\}$ forms a $\delta$-sequence at $t=x$  as $m\to\infty$. In the sequel, the following inequalities that are valid for each $ u \in [0,1] $ will be used: 
\begin{align} \label{eq:4-6}
&| \theta_{m, u}-u | =\frac{|[m u] - m u+1-2 u |}{m+2} \leq  \frac{2}{ m+2}\notag\\
&\sigma^2_{m, u} \leq \frac{u (1-u)}{m+3} \leq \frac{1}{4(m+3)} < \frac{1}{m+2}.
\end{align}
Let us apply  the Taylor series expansion for difference  under the integral in  (\ref{eq:4-4}) and write the left hand side of  (\ref{eq:4-4})  in a  symbolic way as
\begin{align}\label{eq:4-ex}
app(f)-f=I^{}_{10} + I^{}_{01} + I^{}_{20}+ I^{}_{02} + I^{}_{11}.
\end{align}
Now, taking into account (\ref{eq:4-5})-(\ref{eq:4-6}), 
%(\ref{eq:4-4}), 
one can estimate the first four  terms on the right hand side of (\ref{eq:4-ex}) as follows:
\begin{align}\label{eq:4-7}
| I^{}_{10}(x)|&:=|f^{}_{10} (x) \int_{0}^{1} \int_{0}^{1} \, (t - x_1) \beta_{m, x_1}(t)  \beta_{n, x_2}(s)  d t\,  ds|\notag \\
&\leq| f_{10} (x)|\,  | \theta_{m, x_1}-x_1 | \leq  \frac{2\,||f^{}_{10}  ||}{m+2},
\end{align}

\begin{align}\label{eq:4-8}
| I^{}_{01} (x)|&:=|f^{}_{01} (x) \int_{0}^{1} \int_{0}^{1} \, (s - x_2) \beta_{m, x_1}(t)  \beta_{n, x_2}(s)  d t\; ds|\notag \\
&\leq| f^{}_{01} (x)|\,  | \theta_{n, x_2}-x_2 | \leq  \frac{2\,||f^{}_{01}  ||}{n+2},
\end{align}

\begin{align}\label{eq:4-9}
| I^{}_{20} (x)|&:=|\frac{1}{2} \int_{0}^{1} \int_{0}^{1} f^{}_{20} (\tilde t, \tilde s)\, (t - x_1)^2 \beta_{m, x_1}(t)  \beta_{n, x_2}(s)  d t \, ds|\notag \\
&\leq \frac{1}{2}  || f^{}_{20}||\,  \Big ( \sigma^2_{m, x_1} + | \theta_{m, x_1}-x_1 |^2 \Big ) \leq  \frac{\frac{1}{2}\, || f^{}_{20}|| }{m+2} + o\Big (\frac{1}{m}\Big ),
\end{align}

\begin{align}\label{eq:4-10}
| I^{}_{02} (x)|&:=|\frac{1}{2} \int_{0}^{1} \int_{0}^{1} f^{}_{02} (\tilde t, \tilde s)\, (s - x_2)^2 \beta_{m, x_1}(t)  \beta_{n, x_2}(s)  d t\,  ds|\notag \\
&\leq  \frac{1}{2}\, ||f^{}_{02}||\,  \Big ( \sigma^2_{n, x_2} + | \theta_{n, x_2}-x_2|^2 \Big ) \leq  \frac{\frac{1}{2}\, || f^{}_{02}|| }{n+2} + o\Big (\frac{1}{n}\Big )
\end{align}
as $m, n \to\infty$.
To estimate the last term, let us apply  the Cauchy-Schwartz inequality in combination with  (\ref{eq:4-5})-(\ref{eq:4-6}). We obtain
\begin{align}\label{eq:4-11}
&| I^{}_{11} (x)|:=|\frac{1}{2} \int_{0}^{1} \int_{0}^{1} f^{}_{11} (\tilde t, \tilde s)\,(t - x_1)\,  (s - x_2) \beta_{m, x_1}(t)  \beta_{n, x_2}(s)  d t ds|\notag \\
&\leq  \frac{1}{2}\, ||f^{}_{11}||\,  \Big ( \int_{0}^{1}  \, (t- x_1)^2 \beta_{m, x_1}(t) d t \Big )^{1/2}\, \Big (  \int_{0}^{1} \, (s - x_2)^2  \beta_{n, x_2}(s)   ds \Big )^{1/2}\notag \\
&\leq  \frac{1}{2}\, ||f^{}_{11}||\,  \Big ( \frac{1}{m+2}  + \frac{4}{(m+2)^2} \Big )^{1/2} \Big ( \frac{1}{n+2}  + \frac{4}{(n+2)^2} \Big )^{1/2}.
% \frac{1}{2}\, ||f^{\prime\prime}_{11}||\,  | \theta_{n, x}-x | \leq  \frac{2\,||f^{\prime}_{01}  ||}{n+2},
\end{align}
Finally, from inequalities (\ref{eq:4-ex})-(\ref{eq:4-11}) we derive (\ref{eq:4-1}).
\end{proof}
\end{thm}

%%%%%%%%%%%%%%%% convolution f*\phi
In the following statement we use $ \varphi = \varphi_h $ as described in Definition~\ref{defn:molifier}.
Denote the convolution of $f$ and $\varphi$ by $f^{*}:=f \ast \varphi$, while for partial derivatives of $f^{*}$ we will use similar notations as before. For example, we write $f^{*}_{10}:=\partial{f^{*}}/\partial{ x_1}$. 

The following statement about  the  approximation rate of $f$ by $app(f \ast \varphi)$ can be proved.
\begin{thm}\label{conv-moments}
%If  partial derivatives of order two $f_{20},f_{11}, f_{02}$ are  bounded on $ [0,1]^2$,  then
%Let $ f $ be continuously differentiable on $ [0,1]^2 $. 
{Let $ f \in C^2( [0,1]^2) $.} Then
$app( f \ast \varphi ) $ converges uniformly to $ f $  on $ [0,1]^2$ as $ m, n \rightarrow \infty$, and $h\to0$. Furthermore, we have 
\begin{align*}
%\label{eq:4-1}
||app(f \ast \varphi) - f||
& \leq C_1\, h^2+ \frac{ 2\,  || f^{}_{10}|| + \frac{1}{2}\, ||f^{}_{20}|| }{m+2} +  \frac{2\,   ||f^{}_{01}|| + \frac{1}{2}\, ||f^{}_{02}||}{n+2} \notag \\
&+  \frac{\frac{1}{2}\, ||f^{}_{11}||}{\sqrt{(m +2)\,(n+2})}
+ o\Big (\frac{1}{m}\Big ) + o\Big (\frac{1}{n}\Big )+ o\Big (\frac{1}{\sqrt{m\, n}}\Big ). 
\end{align*}
Here $C_1=\frac{\sigma^2}{2} \Big ( ||f^{}_{20}||+ ||f^{}_{02}|| \Big )$ with $\sigma^2=\int t^2\varphi (t) d\,t$. 

In particular, by choosing $m=n$, and $h=\sqrt{\frac{C^{}/C_1}{n+2}}$, we have
\begin{align}
%\label{eq:4-2}
||app(f \ast \varphi) - f||
& \leq  \frac{2\,C^{}}{n+2} + o\Big (\frac{1}{n}\Big )
\end{align}
 as $ n \rightarrow \infty$. 
 %Here  $C^{*}=2 \Big (||f^{*}_{10}|| +  ||f^{*}_{01}||\Big)  + \frac{1}{2}\, \Big ( ||f^{*}_{20}|| + || f^{*}_{11}|| +  ||f^{*}_{02}|| \Big )$.
 %= 2\,  ||f^{\prime}_{10}|| + 2\, ||f^{\prime}_{01}||  + \frac{1}{2}\, \Big ( ||f^{\prime\prime}_{10}|| + || f^{\prime\prime}_{01}|| +  ||f^{\prime\prime}_{11}|| \Big )$.
%\end{thm}

{\rm Note that even though the mollifier $ \varphi $ is smooth, $ f \ast \varphi $ is the convolution of a 2-D function with a 1-D mollifier, which may not be automatically smooth. }

\begin{proof} The moments of $f \ast \varphi$ are related to the moments of 
{ the modified Radon transform (see the first line in Theorem \ref{eq:4-mrt2}), and can be identified  by solving the  corresponding system of equations, similarly to eq. \ref{axb2}.}  
%$f$ and $\varphi$ according to the formula:
%\begin{equation*}
%\label{moment-f-phi}
%\gamma_{j, m}  ( f \ast \varphi)= \sum_{k=0}^j  \sum_{\ell=0}^m\, C(j,k) C(m, \ell)  \,  \gamma_{j-k,  m-\ell} (f)\, \cdot \gamma_{k+ \ell} (\varphi).
%\end{equation*}
{Hence,} application of Theorem~\ref{approx1} provides the  approximation $app(f\ast\varphi)$ of  $f\ast\varphi$ that is based on the moments  $\gamma_{j, m}  ( f \ast \varphi)$. In addition, we have
\begin{align*}
||f-app(f\ast \varphi)||
& \leq ||f- f\ast \varphi ||+ ||f\ast \varphi-app(f\ast\varphi)||.
%& \leq C\, h^2 + o (h^2)+||f\ast \varphi-app(f\ast\varphi)||.
\end{align*}
Now, let us write the difference between $app (f^{*})$ and $f^{*}$ in a similar way as we did in the proof of Theorem~\ref{approx1} for difference between $app (f)$ and $f$. { In particular,} we have:

\begin{align}
app(f^*)-f^*=I^{*}_{10} + I^{*}_{01} + I^{*}_{20}+ I^{*}_{02} + I^{*}_{11}.
\end{align}
%In a similar way we derive the following inequalities for $f\star\phi$:
%Here 
{The rest of the proof mimics the steps used in the proof of  Theorem 5.2.} For example, let us mention that
\begin{align}
| I^{*}_{10}(x)|&:=|f^{*}_{10} (x) \int_{0}^{1} \int_{0}^{1} \, (t - x_1) \beta_{m, x_1}(t)  \beta_{n, x_2}(s)  d t\,  ds|\notag \\
&\leq| f^{*}_{10} (x)|\,  | \theta_{m, x_1}-x_1 | \leq  \frac{2\,||f^{*}_{10}||}{m+2}.
\end{align}
%\begin{align}\label{eq:4-01}
%| I^{*}_{01} (x)|&:=|f^{*}_{01} (x) \int_{0}^{1} \int_{0}^{1} \, (s - x_2)\beta_{m, x_1}(t)  \beta_{n, x_2}(s)  d t\; ds|\notag \\&\leq| f^{*}_{01} (x)|\,  | \theta_{n, x_2}-x_2 | \leq  \frac{2\,||f^{*}_{01}  ||}{n+2},
%\end{align}
%\begin{align}
%\label{eq:4-20}
%| I^{*}_{20} (x)|&:=|\frac{1}{2} \int_{0}^{1} \int_{0}^{1} f^{*}_{20} (\tilde t,\tilde s)\, (t - x_1)^2 \beta_{m, x_1}(t)  \beta_{n, x_2}(s)  d t \, ds|\notag\\&\leq \frac{1}{2}  || f^{*}_{20}||\,  \Big ( \sigma^2_{m, x_1} + |\theta_{m, x_1}-x_1 |^2 \Big ) \leq  \frac{\frac{1}{2}\, || f^{*}_{20}|| }{m+2}+o\Big (\frac{1}{m}\Big ),
%\end{align}
%\begin{align}\label{eq:4-02}
%| I^{*}_{02} (x)|&:=|\frac{1}{2} \int_{0}^{1} \int_{0}^{1} f^{*}_{02} (\tilde t,\tilde s)\, (s - x_2)^2 \beta_{m, x_1}(t)  \beta_{n, x_2}(s)  d t\, ds|\\\notag&\leq  \frac{1}{2}\, ||f^{*}_{02}||\,  \Big ( \sigma^2_{n, x_2} +\theta_{n, x_2}-x_2|^2 \Big ) \leq  \frac{\frac{1}{2}\, || f^{*}_{02}|| }{n+2} +o\Big (\frac{1}{n}\Big ),
%\end{align}
%\begin{align}\label{eq:4-c1}
%&| I^{*}_{11} (x)|:=|\frac{1}{2} \int_{0}^{1} \int_{0}^{1} f^{*}_{11} (\tilde t,\tilde s)\,(t - x_1)\,  (s - x_2) \beta_{m, x_1}(t)  \beta_{n, x_2}(s)  d tds|\notag \\&\leq  \frac{1}{2}\, ||f^{*}_{11}||\,  \Big ( \int_{0}^{1}  \, (t-x_1)^2 \beta_{m, x_1}(t) d t \Big )^{1/2}\, \Big (  \int_{0}^{1} \, (s -x_2)^2  \beta_{n, x_2}(s)   ds \Big )^{1/2}\notag \\&\leq  \frac{1}{2}\,||f^{*}_{11}||\,  \Big ( \frac{1}{m+2}  + \frac{4}{(m+2)^2} \Big )^{1/2} \Big ( \frac{1}{n+2}  + \frac{4}{(n+2)^2} \Big )^{1/2}.
%\end{align}
{Upper bounds similar to (\ref{eq:4-9}) -(\ref{eq:4-11}) can be derived  as well, where instead of  $f_{kj}$ we have $f^{*}_{kj}$ for $k, j =0,1, 2$.} 
Finally, note that 
$||f^{*}_{kj}|| \leq ||f_{kj}||$ for $k+j\leq 2, k,j = 0, 1, 2$. 
%and 
Since $\varphi$ is symmetric and $f$ is smooth, application of the Taylor expansion yields
%we have:
$$
||f- f\ast \varphi || = C_1\, h^2 + o (h^2), 
$$
as $h\to 0$.
\end{proof}
\end{thm}
%%%%%%%%%%%%%%

It is useful to note that Theorem \ref{conv-moments} says that 
we do not need to evaluate the moments of density $ f $; the moments of $ f \ast \varphi $ are derived from the moments of the modified Radon transform using  (\ref{eq:4-mrt2}), 
and this can be used to obtain an approximation of the density function $ f $ by taking  the moments of $ f \ast \varphi_h $ as $ h \to 0 $.  
%Also,  it is easy to note that with additional smoothness assumption on $ f $ in Theorem \ref{conv-moments}, the rate of convergence of the approximation of $f$ can be improved.  

%Also, while the density functions in this section are assumed to be supported in $ [0,1]^2 $, the theorems hold for any moment determinate compactly supported function in $ {\mathbb R}^2 $.  %It is sufficient to have the moments of $ f \ast \varphi_h $ with $ h \to 0 $ to evaluate the moments of $ f $ from the modified Radon transform and reconstruct the density $ f $. 

%%%%%%%%%%%%%%%%%%%%%%%%%%%%%
\section{A Numerical Example}
\label{sec:ANE}
%%%%%%%%%%%%%%%%%%%%%%%%%%%%%

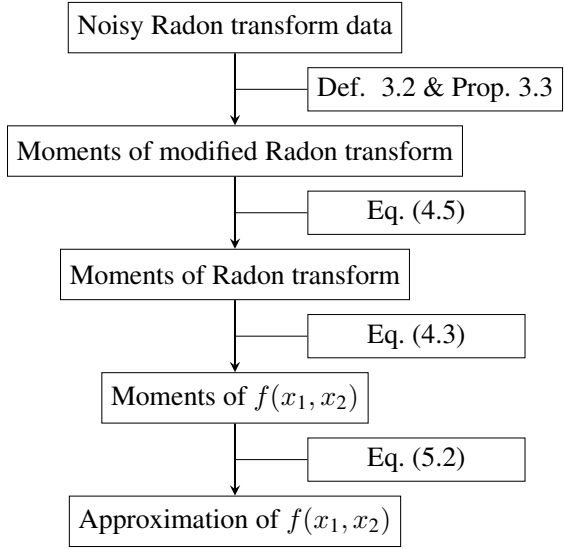
\begin{figure}[tp!]
\begin{tikzpicture}%
  [data/.style=
    {draw,minimum height=0.7cm,minimum width=2cm,align=center},
   filter/.style=
    {draw,minimum height=0.3cm,minimum width=3cm,align=center},
   database/.style=
    {draw,minimum height=0.7cm,minimum width=3cm,align=center},
   flow/.style={thick,-stealth},
   apply/.style={}
  ]
    \node[data] (d0) {Noisy Radon transform data};
  \node[data,below=of d0] (d1) {Moments of modified Radon transform};
  \node[data,below=of d1] (d2) {Moments of Radon transform};
  \node[data,below=of d2] (d3) {Moments of $f(x_1,x_2)$};
  \node[data,below=of d3] (d4) {Approximation of $f(x_1,x_2)$};
    \draw[flow] (d0) -- coordinate(d0d1) (d1);
  \draw[flow] (d1) -- coordinate(d1d2) (d2);
  \draw[flow] (d2) -- coordinate(d2d3) (d3);
  \draw[flow] (d3) -- coordinate(d3d4) (d4);
%\node[filter,right=of d0] at (excl) {Exclusions};
  \node[filter,right=of d0d1] (f1) {Def.~ \ref{defn:MRT} \& Prop.~\ref{modified:form2}};
  \draw[apply] (d0d1) -- (f1);
    \node[filter,right=of d1d2] (f2) {Eq.~\eqref{bbc}};
  \draw[apply] (d1d2) -- (f2);
    \node[filter,right=of d2d3] (f3) {Eq.~\eqref{axb2}};
  \draw[apply] (d2d3) -- (f3);
    \node[filter,right=of d3d4] (f4) {Eq.~\eqref{M-Li}};
  \draw[apply] (d3d4) -- (f4);
  \end{tikzpicture}
  \caption{Outline for our proposed method.}
  \label{fig:fc}
\end{figure}

In this section, we discuss the performance of the proposed procedure to recover a density function from its corresponding modified Radon Transform moments. We give a bird's eye view of the inversion algorithm as a flow chart (see Figure~\ref{fig:fc})
that describes the different steps in the reconstruction.

To provide an accurate simulation and minimize the pollution inherent in floating-point calculation with limited precision, we developed a computer code that heavily utilizes the GNU MPFR Library (https://www.mpfr.org). This is a C library for multiple-precision floating-point computations with correct rounding. The linear algebraic equations in the proposed procedure are solved using solvers in Eigen, which is  a high-level C++ library of template headers for linear algebra, matrix and vector operations, geometrical transformations, numerical solvers and related algorithms (see http://eigen.tuxfamily.org). Application of Eigen is made possible by an MPFR C++ wrapper (see http://www.holoborodko.com/pavel/mpfr/).

\begin{figure}[H]
\centering
\includegraphics[scale=0.65]{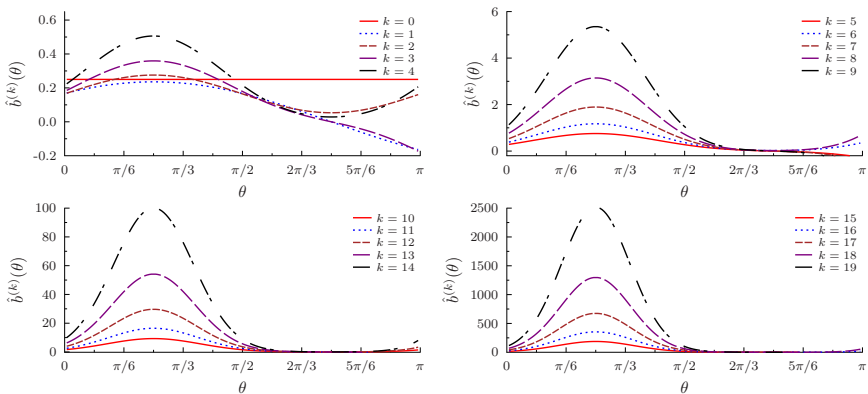}
\caption{A sample of moments of the modified Radon Transform, $\hat{b}^{(k)}(\theta)$.} \label{fig:modrdm01}
\end{figure}

\begin{figure}[H]
\centering
\includegraphics[scale=0.65]{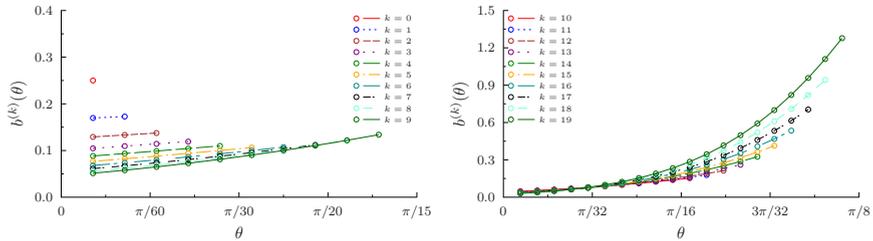}
\caption{The resulting moments of Radon Transform, $b^{(k)}(\theta)$  obtained from the data depicted in Figure~\ref{fig:modrdm01}.} \label{fig:rdm01}
\end{figure}

\begin{figure}[H]
\centering
\includegraphics[scale=0.8]{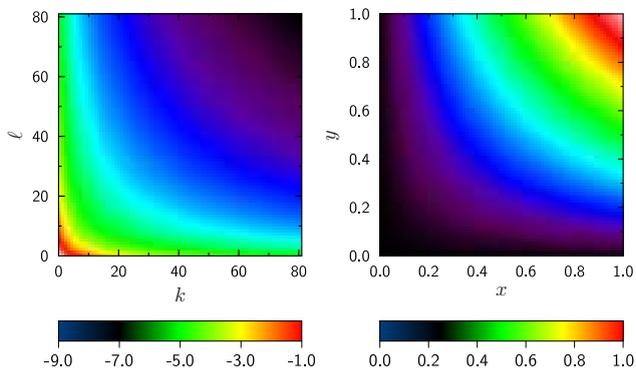}
\caption{Left: Logarithmic plot of density moments, $\gamma_{k,\ell}(f)$, right: the corresponding predicted density $f(x,y)$.} \label{fig:den01}
\end{figure}

The target density function is a known function $f(x,y) = xy$. To illustrate the procedure, we assume availability of a set of moments of the modified Radon transform data,  $\hat{b}^{(k)}(\theta)$, for $\theta \in (0,\pi/4)$, $\theta \in (\pi/4,\pi/2)$, $\theta \in (\pi/2,3\pi/4)$, and $\theta \in (3\pi/4,\pi)$, each of them contains $41$ discrete points (a total of $164$ points). A sample of this data is depicted in Figure~\ref{fig:modrdm01}. Equation \eqref{bbc} is then inverted to obtain $b^{(k)}(\theta)$, whose result is shown in Figure~\ref{fig:rdm01}. It is clear that as the order of moments increases, the magnitude of these moments increases as well. At the next stage, a series of inversions of Equation \eqref{axb2} was done to obtain the moments  $\gamma_{k,\ell}(f)$ of function $ f $ and in turn Equation~\eqref{M-Li} was used to get the target density approximation. The results of these calculations are depicted in Figure~\ref{fig:den01}.

\begin{figure}[H]
\centering
\includegraphics[scale=1.]{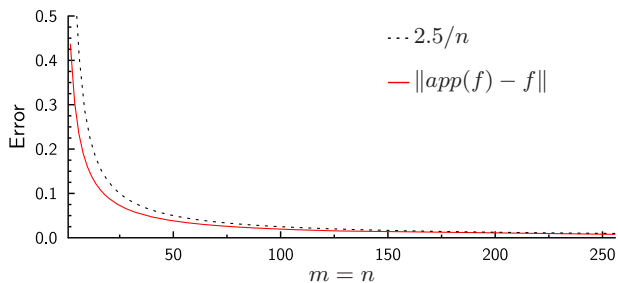}
\caption{Convergence of the density approximation $||app (f) - f ||$ versus number of moments used.} \label{fig:conden}
\end{figure}

Convergence behavior of the density approximation with respect to the moments order $m$ and $n$ is shown in Figure~\ref{fig:conden}. A bound of this error in the form of $\mathcal{O}(1/n)$ is also plotted in this figure. This figure confirms the theoretical finding established in Theorem~\ref{approx1}.


\section{Concluding Remarks}
\label{sec:concl}

While there are many excellent monographs and papers on the Radon transform and its applications to tomography,  because of the significance of this transform, new methods are being continuously developed. A common theme among these methods are:
\begin{enumerate}
    \item How can the reconstruction be made specific to highlight specific features in the image?
    \item How can noise or other artifacts be suppressed?
    \item How can the reconstruction be performed optimally from fewer projections?
\end{enumerate}
This paper deals with all of these issues at a theoretical level.  By moving from Fourier methods underlying the standard FBP and ART algorithms to using moment methods, we show how mollification of the Radon transform is transported into the moment problem, and derive explicit relationships between the moments of Radon transform, moments of its mollified version, and moment of the original density function. We also show how these reconstructions from the moments of the modified transform converge uniformly to the original density function (and not just the mollified density function). A numerical example provides details of this approximation, and verifies the accuracy of the theoretically derived algorithm and its convergence rate.

\noindent We have left the extensive study of finding optimal mollifiers for individual applications, numerical results of density patterns with discontinuities and/or more anthropomorphically realistic patterns, and generalization of these methods to a succession of future papers.

%%%%%%%%%%%%%%%%%%%%%%%%%%%%%%%%%%%%%%%%%%%%%%%%%%%%%%%% %%%%%%%%%  Reference
%%%%%%%%%%%%%%%%%%%%%%%%%%%%%%%%%%%%%%%%%%%%%%%%%%%%%%%%

\bibliographystyle{plain}
\bibliography{refs}

\end{document}